# Distributional Robust Portfolio Construction based on Investor Aversion

Xin ZHANG[a]

[a]School of Management and Economics, BeiHang University, Beijing 100083, PR China

***Abstract:*** *In behavioral finance, aversion affects investors' judgment of future uncertainty when profit and loss occur. Considering investors' aversion to loss and risk, and the ambiguous uncertainty characterizing asset returns, we construct a distributional robust portfolio model (DRP) under the condition that the distribution of risky asset returns is unknown. Specifically, our objective is to find an optimal portfolio of assets that maximizes the worst-case utility level on the Wasserstein ball, which is centered on the empirical distribution of sample returns and the radius of the ball quantifies the investor's ambiguity level. The model is also formulated as a mixed-integer quadratic programming problem with cardinality constraints. In addition, we propose a hybrid algorithm to improve the efficiency of the solution and make it more suitable for large-scale problems. The performance of the DRP model since the outbreak of COVID-19 is empirically investigated, verifying that the DRP model is effective in avoiding systemic risk and achieving excess returns over general asset allocation strategies during periods of sharp financial market volatility.*

## 0 Introduction:

The investment portfolio is a collection of stocks, bonds, financial derivatives, etc. held by the investor or financial institution. The purpose is to diversify risk. One important aspect is to consider how to combine risky assets. Since any combination of two poorly or negatively correlated assets will yield a riskier return than the risky return of the assets alone, constantly combining poorly correlated assets can keep the portfolio's efficient frontier away from risk. Since mean-variance model (Markowitz 1952), which pioneered the analysis of financial mathematics, the theory of portfolio optimization has evolved in terms of models and practical applications, but still faces many challenges. Tversky (1979) proposed the prospect theory, and since then people have started to pay attention to the impact of aversion on portfolio selection. Dekel (1989) examine the relationship between risk aversion and portfolio diversification when preferences over probability distributions of wealth do not have an expected utility representation. Ramaswami et al. (1992) examine the relationship between asset holdings and portfolio objectives. Jarrow and Zhao (2006) compare mean-variance (M-V) and mean-lower partial moment (M-LPM) optimal portfolios under nonnormal asset return distributions. Using household survey data (Dimmock and Kouwenberg 2010) obtain direct measures of each surveyed household's loss-aversion coefficient from questions involving hypothetical payoffs. Warnick et al. (2011) argue that ambiguity aversion is just as relevant to their decision-making process because they are uncertain about the yield distributions generated by new technologies. Hanna et al. (2011) discuss the relationship between risk aversion and portfolio recommendations based on an expected utility approach, review selected empirical research on risk tolerance, and propose to separate risk capacity, expectations, and other factors from the concept of risk tolerance. Baltzer, Stolper, and Walter (2015) analyze the effect of geographic proximity on individual investors' portfolio choice. Sun et al. (2015) introduce a new portfolio selection method, the portfolio selection problem is formulated as a bi-criteria optimization problem that maximizes the expected return of the portfolio and minimizes the maximum individual risk of the assets in the portfolio.

In general, when describing expected returns, scholars usually assume that the expected return of an asset is certainty or follows a certain deterministic distribution. Füllbrunn et al. (2014) point that investors are ambiguity averse to return estimation errors if they cannot estimate the true value of the risky asset's return with certainty. Alonso, Prado, and Fluctuations (2008) also argue that investors are ambiguity averse to uncertainty and that this characteristic is prevalent in financial markets. Many scholars use the idea of robust optimization in describing ambiguity aversion characteristics. For example, Garlappi, Uppal, and Wang (2007) use the assumption that the mean of returns belongs to the ellipsoidal uncertainty set to represent investors' ambiguity aversion and study the mean-variance robust optimization problem based on ambiguity aversion characteristics in terms of parameters and model uncertainty, respectively. On this basis, Pınar (2014) and Fabretti, Herzel, and Pınar (2014) study the equilibrium market pricing under ambiguity aversion mean-variance conditions and the delegated portfolio selection problem for ambiguity averse investors, respectively. Chen, Ju, and Miao (2014) used a recursive ambiguity utility model to describe the ambiguity aversion characteristics of investors and pointed out that the robust investment strategy of ambiguity aversion investors is more conservative than the Bayesian strategy. However, all the above robust optimization ideas assume that the asset return probability distribution is known or follows a normal distribution when describing ambiguity aversion, and do not consider investors' ambiguity aversion to uncertainty of asset return distribution. Ellsberg (1961) and Liu and Control (2011) both point out that the case of unknown probability distribution of asset returns should be considered when describing the characteristics of ambiguity aversion. Therefore, borrowing from the idea of distribution robustness which reflects the uncertainty of asset return distribution (Delage and Ye 2010; Goh and Sim 2010; Long and Qi 2014), it is more consistent with the real psychology of investors to introduce the assumption of unknown return distribution in robust optimization to describe investors' ambiguity aversion characteristics.

However, traditional robust optimization is too conservative and always considers the worst-case scenario (Fabretti, Herzel, and Pınar 2014; Natarajan, Pachamanova, and Sim 2008; Gregory, Darby-Dowman, and Mitra 2011; Chen, He, and Zhang 2011; Sehgal and Mehra 2020). In portfolio building blocks, higher risk also implies higher return, so the overly conservative strategy is not suitable for general situations. Stochastic optimization strategies often ignore the properties of the return distribution, resulting in unconvincing decisions (Kraft 2004; Kamil, Mustafa, and Ibrahim 2009; Cui et al. 2020; Guigues 2021). The distributional robust optimization method based on Wasserstein ambiguity ensemble combines the advantages of traditional robust optimization and stochastic optimization by controlling the degree of ambiguity through the radius of ambiguity sphere, which can consider the optimal solution of the model under the worst distribution (Gao and Kleywegt 2016; Hanasusanto and Kuhn 2018; Mohajerin Esfahani and Kuhn 2015). Researchers have applied it to the portfolio field, which can well describe the stochasticity and ambiguity of returns, and at the same time better quantify the ambiguity aversion characteristics of investors. Pflug, Pohl, and Analysis (2018) show that in this case portfolio concentration becomes optimal as the uncertainty with respect to the estimated dependence structure increases. Therefore, Roveto, Mieth, and Dvorkin (2020) develop a means of co-optimizing the value-at-risk ($VaR$) level associated with the $CVaR$ to guarantee resilience in probable cases while providing a measure of the average violation in less probable cases. Since optimal portfolio strategy depends heavily on the distribution of uncertain returns, Du, Liu, and Liu (2020) propose a new method for the portfolio optimization problem with respect to distribution uncertainty. Blanchet, Chen, and Zhou (2021) revisit Markowitz's mean-variance portfolio selection model by considering a distributional robust version, where the region of distributional uncertainty is around the empirical measure and the discrepancy between probability measures is dictated by the so-called Wasserstein distance. Pesenti and Jaimungal (2021) study the problem of active portfolio management where an investor aims to outperform a benchmark strategy's risk profile while not deviating too far from it. Campbell and Wong (2021) develop a concrete and fully implementable approach to the optimization of functionally generated portfolios in stochastic portfolio theory. The focus is on the quality of solutions achieved as determined by the Normalized and Complementary Wasserstein Distance $\eta$ , which Baker and Radha (2022) present in a manner to expose the QAOA as a transporter of probability.

In this paper, we propose a novel data-driven Distributional Robust Portfolio（DRP）model with a metric-based ambiguity set for the unknown distribution of asset return rate, which is defined as a Wasserstein ball centered at the empirical distribution over the finite sample dataset and the radius of the ball reflects the investor's ambiguity perception of the return rate distribution. That is, the higher the investor's ambiguity to returns, the lower the credibility of the empirical distribution and the larger the radius of the ball; conversely, the lower the ambiguity aversion, the smaller the radius of the ball. The objective of the portfolio is to maximize a utility function incorporating loss and risk aversion coefficients, the portfolio does not allow short selling and is subject to cardinality constraints, and the final portfolio problem is modeled as a mixed integer quadratic programming. The sections of this paper are organized as follows: Section 1, we introduce the concept of investors' main aversions and construct utility functions based on loss and risk aversions. Section 2, a distributional robust optimization portfolio model (DRP) is constructed based on Wasserstein ambiguity sets. In Section 3, by considering the Dual Theory, the model is reformulated as a solvable mixed-integer quadratic programming, in addition to giving the worst-case distribution. Section 4, problem decomposition and optimality conditions for sub-quadratic programming problems Section 5, we propose a hybrid algorithm to handle the DRP and improve the computational efficiency. Section 6 is a sensitivity and model robustness test for the aversion coefficient. Section 7, where we set the rolling time window based on the historical return of SP&500 during Corona Virus Disease 2019 (OVID-19) and compare the performance of DRP with common asset allocation strategies to verify the feasibility of DRP strategies. Section 8 is some concluding work and outlook.

## 1 Investor Emotion and Utility Functions:

Behavioral finance is the integration of theories from psychology (especially behavioral science) into finance, which studies and predicts the development of financial markets in terms of micro-individual behavior and the psychological and other motivations that produce such behavior. In the field of investment portfolio, it is often considered that investors' aversion can lead to irrational decisions. Common investor aversions include:

($i$). Risk aversion:

In economics and finance, risk aversion is the tendency of people to prefer outcomes with low uncertainty to those outcomes with high uncertainty, even if the average outcome of the latter is equal to or higher in monetary value than the more certain outcome. Risk aversion explains the inclination to agree to a situation with a more predictable, but possibly lower payoff, rather than another situation with a highly unpredictable, but possibly higher payoff (Kahneman and Tversky 2013).

($ii$). Loss aversion:

Loss aversion is the tendency to prefer avoiding losses to acquiring equivalent gains. The principle is prominent in the domain of economics. What distinguishes loss aversion from risk aversion is that the utility of a monetary payoff depends on what was previously experienced or was expected to happen.

($iii$). Ambiguity aversion:

In decision theory and economics, ambiguity aversion (also known as uncertainty aversion) is a preference for known risks over unknown risks (Gollier 2005). An ambiguity-averse individual would rather choose an alternative where the probability distribution of the outcomes is known over one where the probabilities are unknown.

Incorporating aversion into investor utility, we construct the utility function:

$$\begin{cases} g(R_p) = R_p(\tilde{\xi}) - \varphi\left[\hat{R} - R_p(\tilde{\xi})\right]^+ - \dfrac{1}{2}\mathcal{A}\sigma_p \\ \tilde{\xi} \in \mathcal{F}(\theta) \end{cases}$$

Where $R_p$ is the portfolio return, $\varphi \geq 0$ is the loss aversion coefficient, $\hat{R}$ is the reference point for defining losses and gains, $\left[\hat{R} - R_p\right]^+$ is defined as the greater of $\hat{R} - R_p$ and 0, which is $\max\{\hat{R} - R_p, 0\}$. The function is linear in both gain part and loss part. $\mathcal{A}$ is the coefficient to measure the level of risk aversion. In addition, to reflect investors' aversion to ambiguity returns, returns are defined within an ambiguity set $\mathcal{F}(\theta)$, the size of which is controlled by the ambiguity aversion coefficient $\theta$.

## 2 Distributional Robust Portfolio:

Distributional robust optimization is a paradigm for decision making under uncertainty where the uncertain problem data are governed by a probability distribution that is itself subject to uncertainty. The core of distributional robust optimization is in the construction of data-driven ambiguity set, where the random variables can be discrete or continuous. For the continuous case, familiar ambiguity set are constructed in two ways: ($i$) Constructed based on moment information: It contains all distributions that satisfy the moment constraints (moments of each order) of the random variables; ($ii$) Metric-based ambiguity set: Using probabilistic distance functions (e.g., Prokhorov metric, K-L divergence, Wasserstein metric etc.) to define the ambiguity set as a probability space distribution sphere, where the center of the sphere is the empirical distribution of the historical data, and then the sphere with theta as the ambiguity radius contains the unknown true distribution of the random variables.

### 2.1 Wasserstein Ambiguity Set

We consider a data-driven setting as in Kuhn et al. (2019) on the design of a Wasserstein ambiguity set centered around the empirical distribution $\mathbb{P}_e(\xi) = \frac{1}{s}\sum_{i=1}^{s} \delta_{\hat{\xi}_i}$, where $\delta$ is the Dirac measure (i.e. $\delta_{\hat{\xi}_i}(A) = 1$, if $\hat{\xi}_i \in A$, and otherwise $\hat{\xi}_i \notin A$, $\delta_{\hat{\xi}_i}(A) = 0$, for any set $A \in \sigma(\Omega)$). When the sample is large enough the true distribution $\mathbb{P}$ will converge to the empirical distribution $\mathbb{P}_e$.

Given a tractable distance metric $d_w: \mathcal{M}(\Xi) \times \mathcal{M}(\Xi) \mapsto \mathfrak{R}^+$, the Wasserstein metric (a.k.a Kantorovich-Rubinstein metric) between any two distributions $\mathbb{P}_1$ and $\mathbb{P}_2$ is defined via an optimization problem. And for $p \in [1,\infty)$, the type-p Wasserstein metric between two distributions $\mathbb{P}$ and $\widehat{\mathbb{P}}$ for a given distance metric $\rho$ is defined as:

$$d_w(\mathbb{P}_1,\mathbb{P}_2) \triangleq \inf\left\{\int_{\Xi\times\Xi} \rho^p(\xi_1,\xi_2)\mathcal{K}(d\xi_1,d\xi_2):\ \int_{\Xi} \mathcal{K}(\xi_1,d\xi_2) = \mathbb{P}_1(\xi_1), \int_{\Xi} \mathcal{K}(d\xi_1,\xi_2) = \mathbb{P}_2(\xi_2) \right\}$$

where $\xi_1 \sim \mathbb{P}_1$, $\xi_2 \sim \mathbb{P}_2$, $(\Xi,\rho)$ is a Polish metric space, $K:\Xi \times \Xi \mapsto \mathfrak{R}^+$ is the joint distribution of $\mathbb{P}_1 \in \mathcal{M}(\Xi)$ and $\mathbb{P}_2 \in \mathcal{M}(\Xi)$. Moreover, $\rho^p(\xi_1,\xi_2) = \|\xi_1 - \xi_2\|_p$ where $\|\bullet\|_p$ represents $l_p - norm$ on $\mathfrak{R}^n$. We use type-1 Wasserstein metric to construct ambiguity sets.

The type-p Wasserstein ambiguity set is then defined by

$$\mathcal{F}_W(\theta) = \left\{\widehat{\mathbb{P}} \in \mathcal{M}(\Xi) \;\middle|\; \begin{array}{c} \xi \sim \widehat{\mathbb{P}} \\ d_w(\mathbb{P}_e,\widehat{\mathbb{P}}) \le \theta \end{array} \right\}$$

which is a ball of radius $\theta \ge 0$, $\mathbb{P}_e$ is the center of the ball, $\widehat{\mathbb{P}}$ is an estimate of the true distribution $\mathbb{P}$.

2.2. DRP Model

Considering the perspective of an investor with multiple aversion characteristics, whose aversion coefficient $\mathcal{A} > 0$ and loss aversion coefficient $\varphi > 0$. $\hat{R} \in \mathfrak{R}$ is a given reference point for measuring loss and gain, and $\xi$ is the stochastic rate of return. $\hat{\xi}_i$ is the $i$ historical sample of the $\xi$. The number of historical data samples is $[S]$ and the number of total assets is $[N]$. Then, the distributional robust portfolio model (DRP) can be modeled as:

$$\underset{x}{maximize} \inf_{\xi\sim P,\ P\in \mathcal{F}_W(\theta)} \left\{\mathbb{E}_{\mathbb{P}}\left[\xi^T x - \varphi\left[\hat{R} - \xi^T x\right]^+\right]\right\} - \frac{1}{2}\mathcal{A}\, x^T\hat{\Sigma}x \qquad (a)$$

$$s.t. \quad \sum_{j=1}^{n} x_j = 1$$

$$\sum_{j=1}^{n} y_j = k \qquad (Cardinality\ Constraint)$$

$$x_j \le y_j, \forall j \in [N]$$

$$x_j \in [0,1], \forall j \in [N] \qquad (Short\ selling\ is\ not\ allowed)$$

$$y_j \in \{0,1\}, \forall j \in [N]$$

**Parameters:**

·S: The number of historical samples of asset returns, for $\forall j \in [N]$ there are $s$ samples with a sample index of $\forall i \in [S]$.

·N: Number of total assets, with assets indexed as $\forall j \in [N]$.

·$\hat{\xi}_i$,$\hat{\xi}_{ij}$: $\hat{\xi}_i$ is the rate of return of assets in the sample group $i$, $\hat{\xi}_{ij}$ is the rate of return on asset $j$ in sample group $i$ ,$\forall j \in [N]$

·$\xi$: The vector of stochastic returns on assets.

·$\hat{\Sigma}$: Covariance matrix between assets, which is estimated based on historical data samples.

$\cdot\varphi$: The loss aversion coefficient of the investor, for the loss aversion $\varphi \geq 0$.
$\cdot\mathcal{A}$: The risk aversion coefficient of the investor, for the risk aversion $\mathcal{A} > 0$.
$\cdot\theta$: The radius of the Wasserstein ambiguity set, which controls the size of the ambiguity set, i.e., the Wasserstein distance between two distributions is not greater than $\theta$. In addition, $\theta$ reflects the investor's ambiguity aversion, the larger $\theta$, the higher the fuzziness of the return.
$\cdot\widehat{R}$: $\widehat{R} \in \mathfrak{R}$ is a given reference point for measuring loss and gain, it can be artificially set, for example, based on a comparison of the past wealth with the current wealth; or it can be a criterion for distinguishing between losses and gains with a certain objective in mind. Here we use the index return as a reference point
$\cdot k$: Cardinality constraints for limiting portfolio size.
Decision Variables:
$\cdot\ x$: Continuous decision variables, $x = [x_1,x_2,...,x_n]^T$ is the weight vector of investment decision makers.
$\cdot y$: 0-1 logical decision variables, $y = [y_1,y_2,...,y_n]^T$ and $y_j = 1$ means asset $j$ is selected for the portfolio.

The objective function of the model is to maximize the investor's sentiment utility under the worst distribution, and the constraints mainly include the cardinality constraint and the position constraint.

## 3 Reformulating the DRP model and Worst-case Distribution Derivation:

3.1 Reformulating the DRP model:

In order to assess the maximum emotional utility of the investment decision maker, we first need to consider the worst-case distribution $\mathbb{P}^*$ part in $(a)$:

$$\inf_{\xi\sim P,\ P\in\mathcal{F}_W(\theta)} \left\{\mathbb{E}_{\mathbb{P}}\left[\xi^T x - \varphi\left[\widehat{R} - \xi^T x\right]^+\right]\right\}$$

To facilitate writing, make $h(x,\varphi,\widehat{R},\xi) = \xi^T x - \varphi\left[\widehat{R} - \xi^T x\right]^+$

Then,

$$\mathbb{E}_{\mathbb{P}^*}\left[h(x,\varphi,\widehat{R},\xi)\right] = \inf_{\mathbb{P}\in\mathcal{F}_W(\theta)} \mathbb{E}_{\mathbb{P}}\left[h(x,\varphi,\widehat{R},\xi)\right]$$

It is easy to see that $\inf_{\mathbb{P}\in\mathcal{F}_W(\theta)} \mathbb{E}_{\mathbb{P}}\left[h(x,\varphi,\widehat{R},\xi)\right]$ can be derived by solving the following cone programming:

$$\underset{K(\xi,\hat{\xi}_i)\geq 0}{minimize} \int_{\Xi} \sum_{i=1}^{S} h(x,\varphi,\widehat{R},\xi)\mathcal{K}(d\xi,\hat{\xi}_i) \qquad \#(b)$$

$$s.t. \qquad \int_{\Xi} \mathcal{K}(d\xi,\hat{\xi}_i) = \frac{1}{S}, \forall i \in [S]$$

$$\int_{\Xi} \sum_{i=1}^{S} \rho^p(\xi,\hat{\xi}_i)\mathcal{K}(d\xi,\hat{\xi}_i) \leq \theta$$

We introduce auxiliary variables $\lambda$,$\nu$ and the Lagrange function.

$$\mathcal{L}(\xi,\lambda,\nu) = \int_{\Xi} \sum_{i=1}^{S} h(x,\varphi,\widehat{R},\xi)\mathcal{K}(d\xi,\hat{\xi}_i) + \int_{\Xi} \sum_{i=1}^{S} \nu_i\mathcal{K}(d\xi,\hat{\xi}_i) + \int_{\Xi} \sum_{i=1}^{S} \lambda\rho^p(\xi,\hat{\xi}_i)\mathcal{K}(d\xi,\hat{\xi}_i) - \frac{1}{S}\sum_{i=1}^{S} \nu_i - \lambda\theta$$

The original problem $(b)$ is equivalent to $\min_{\xi\in\Xi}\max_{\lambda,\nu} \mathcal{L}(\xi,\lambda,\nu)$

It follows that the Lagrange dual problem can be represented as:

$$\max_{\lambda,\nu}\min_{\xi\in\Xi} \mathcal{L}(\xi,\lambda,\nu) = \max_{\lambda,\nu}\min_{\xi\in\Xi} \int_{\Xi} \sum_{i=1}^{S} \left(h(x,\varphi,\widehat{R},\xi) + \nu_i + \lambda\rho^p(\xi,\hat{\xi}_i)\right)\mathcal{K}(d\xi,\hat{\xi}_i) - \frac{1}{S}\sum_{i=1}^{S} \nu_i - \lambda\theta$$

Consequently, the dual problem is given as:

$$\underset{\lambda,\nu}{maximize} - \frac{1}{s}\sum_{i=1}^{s} \nu_i - \lambda\theta \qquad \#(c)$$

$$s.t. \qquad - h\left(x,\varphi,\widehat{R},\xi\right) - \lambda\rho^p\left(\xi,\hat{\xi}_i\right) \leq \nu_i, \forall \xi \in \Xi, i \in [S]$$

$$\lambda \geq 0$$

Consider the $(b)$ and its dual problem $(c)$. If $\theta > 0$, there exists a strictly feasible solution $\mathcal{K} = \mathbb{P}_e \times \mathbb{P}_e$ to the $(b)$. Thus, the Slater condition for their strong duality holds (Shapiro, 2001). If $\theta = 0$, the Wasserstein ball $\mathcal{F}_W(\theta)$ reduces to a singleton $\{\mathbb{P}_e\}$ and $(b)$ changes to a sample average problem $\frac{1}{S}\sum_{i=1}^{S} h\left(x,\varphi,\widehat{R},\xi\right)$ Indeed, any feasible solution to the dual problem $(c)$ satisfies that $\nu_i \geq - h\left(x,\varphi,\widehat{R},\xi\right)$ with $\xi = \hat{\xi}_i$ and $\nu_i \geq 0$ with $\xi \neq \hat{\xi}_i$ when $\lambda$ tends to infinity.

Accordingly, the optimal value of problem $(c)$ reduces to the sample average problem $\frac{1}{s}\sum_{i=1}^{s} h\left(x,\varphi,\widehat{R},\xi\right)$ as well. Overall, there is no duality gap between $(b)$ and $(c)$ under any case. Thus, it is sufficient to solve the dual problem $(c)$.

To further simplify, we introduce Lemma 1

***Lemma1:*** *For any* $w \in \Re^n$, *it holds that*

$$\sup_{x\in \Re^n}\left\{w^T x - \lambda \parallel x\parallel_p\right\} = \sup_{x\in \Re^n}\left\{\left(\parallel w\parallel_q - \lambda\right) \parallel x\parallel_p\right\},$$

*Where* $\parallel \bullet \parallel_q$ *is the dual of* $l_p - norm$, *i.e.,1/p + 1/q = 1.*

***Proof 1:***

$$\sup_{x\in \Re^n}\left\{w^T x - \lambda \parallel x\parallel_p\right\} = \sup_{t\geq 0}\sup_{\parallel x\parallel_p = t}\left\{w^T x - \lambda \parallel x\parallel_p\right\}$$

$$= \sup_{t\geq 0}\sup_{\parallel x\parallel_p = t}\left\{w^T x - \lambda t\right\} = \sup_{t\geq 0}\left\{t \parallel w\parallel_q - \lambda t\right\}$$

$$= \sup_{\parallel x\parallel_p \geq 0}\left\{\left(\parallel w\parallel_q - \lambda\right) \parallel x\parallel_p\right\}$$

$$= \sup_{x\in \Re^n}\left\{\left(\parallel w\parallel_q - \lambda\right) \parallel x\parallel_p\right\}$$

Since $h\left(x,\varphi,\widehat{R},\xi\right) = \xi^T x - \varphi\left[\widehat{R} - \xi^T x\right]^+$ and $\rho^1\left(\xi,\hat{\xi}_i\right) = \left\|\xi - \hat{\xi}_i\right\|_1$, the constraint in $(c)$ amounts to:

$$\begin{cases} \sup_{\xi\in\Xi}\left\{\varphi\left(\widehat{R} - \xi^T x\right) - \xi^T x - \lambda\left\|\xi - \hat{\xi}_i\right\|_1\right\} \leq \nu_i, \widehat{R} - \xi^T x \geq 0 & (1d) \ \#\#\#\#\#\# \\ \sup_{\xi\in\Xi}\left\{- \xi^T x - \lambda\left\|\xi - \hat{\xi}_i\right\|_1\right\} \leq \nu_i, \qquad \widehat{R} - \xi^T x \leq 0 & (2d) \end{cases}$$

For inequality $(1d)$: we denote $\Delta u_i = \xi - \hat{\xi}_i$, and re-express the left-hand side as:

$$\sup_{\xi\in\Xi}\left\{\varphi(\widehat{R} - x^T\left(\Delta u_i + \hat{\xi}_i\right) - x^T\left(\Delta u_i + \hat{\xi}_i\right) - \lambda\|\Delta u_i\|_1\right\}$$

$$= \sup_{\xi\in\Xi}\left\{\varphi\widehat{R} - \varphi x^T\Delta u_i - \varphi x^T\hat{\xi}_i - x^T\Delta u_i - x^T\hat{\xi}_i - \lambda\|\Delta u_i\|_1\right\}$$

$$= \sup_{\xi\in\Xi}\left\{- (1 + \varphi)x^T\Delta u_i - \lambda\|\Delta u_i\|_1\right\} + \varphi\widehat{R} - (1 + \varphi)x^T\hat{\xi}_i$$

$$= \sup_{\xi\in\Xi}\left\{[(1 + \varphi)\|x\|_\infty - \lambda] \cdot \|\Delta u_i\|_1\right\} + \varphi\widehat{R} - (1 + \varphi)x^T\hat{\xi}_i$$

$$= \begin{cases} \varphi\widehat{R} - (1 + \varphi)x^T\hat{\xi}_i, & (1 + \varphi)\|x\|_\infty \leq \lambda \\ \infty, & (1 + \varphi)\|x\|_\infty > \lambda \end{cases}$$

Similarly, for inequality $(2d)$, we obtain:

$$\sup_{\xi\in\Xi}\left\{\xi^T x - \lambda\left\|\xi - \hat{\xi}_i\right\|_1\right\} = \begin{cases} - x^T\hat{\xi}_i, & \|x\|_\infty \leq \lambda \\ \infty, & \|x\|_\infty > \lambda \end{cases}$$

So, the dual problem $(c)$ can be reformulated as:

$$\underset{\lambda,\nu}{maximaze} - \frac{1}{S}\sum_{i=1}^{S} \nu_i - \lambda\theta \quad \#(e)$$

$$\begin{aligned} s.t. \quad & \nu_i + (1+\varphi)\hat{\xi}_i^T x - \varphi\hat{R} \geq 0, i \in [S] \\ & \nu_i - \hat{\xi}_i^T x \geq 0, i \in [S] \\ & \|x\|_\infty \leq \frac{\lambda}{1+\varphi} \\ & \lambda \geq 0 \end{aligned}$$

Then, DRP can be reformulated as:

$$\underset{x,y,\lambda,\nu}{maximize} - \frac{1}{S}\sum_{i=1}^{S} \nu_i - \lambda\theta - \frac{1}{2}\mathcal{A}\, x^T\hat{\Sigma}x \quad (f)$$

$$\begin{aligned} s.t. \quad & \sum_{j=1}^{n} x_j = 1 \\ & \sum_{j=1}^{n} y_j = k \quad (Cardinality\ Constraint) \\ & x_j \leq y_j, \forall j \in [N] \\ & \nu_i + (1+\varphi)\hat{\xi}_i^T x - \varphi\hat{R} \geq 0, \forall i \in [S] \\ & \nu_i - \hat{\xi}_i^T x \geq 0, \forall i \in [S] \\ & \|x\|_\infty \leq \frac{\lambda}{1+\varphi} \quad (Equivalent\ to: -\frac{\lambda}{1+\varphi} \leq x_j \leq \frac{\lambda}{1+\varphi}, \forall j \in [N]) \\ & \lambda \geq 0 \\ & x_j \in [0,1], \forall j \in [N] \quad (Short\ selling\ is\ not\ allowed) \\ & y_j \in \{0,1\}, \forall j \in [N] \end{aligned}$$

Then, the infinite dimensional programing problem is reformulated as a mixed integer quadratic convex optimization problem which can be solved by the solver. The problem is similar to a knapsack problem, which is a typical NP-hard problem, so the exact solution cannot be found in polynomial time.

3.2 The worst-case distribution of the DRP

We can derive the worst-case distribution of stochastic returns for any given loss aversion characteristic $\{x,y,\varphi,\hat{R}\}$, $x,\varphi,\hat{R} \in R$.

Based on the study (Wang et al. 2020), it can be proved that：

***Lemma2:***

*For the investor with deterministic aversion characteristics, if given a set of feasible solutions* $\{x,y|x_j \in [0,1], y_j \in \{0,1\}, \forall j \in [N]\}$,

$\inf_{\mathbb{P}\in \mathcal{F}_W(\theta)} \mathbb{E}_{\mathbb{P}}[h(x,\varphi,\hat{R},\xi)] \equiv \inf_{\xi \in \mathcal{B}} \frac{1}{S}\sum_{i=1}^{S} h(x,\varphi,\hat{R},\xi_{(i)})$ *always holds.*

*Where* $\mathcal{B} = \{\xi_{(1)},\xi_{(2)},\ldots,\xi_{(i)} | \frac{1}{S}\sum_{i=1}^{S} \rho^1(\xi_{(i)},\hat{\xi}_i) \leq \theta, \xi_{(i)} \in \Xi\}$.

***Lemma3:***

*Given any feasible pair of* $\{x,y,\varphi,\hat{R}\}$, *let* $\xi^{\{x,y,\varphi,\hat{R}\}} = (\xi_{(1)}^{\{x,y,\varphi,\hat{R}\}},\xi_{(2)}^{\{x,y,\varphi,\hat{R}\}},\ldots,\xi_{(S)}^{\{x,y,\varphi,\hat{R}\}})$ *be an optimal solution of the optimization problem. Then the following distribution:*

$\mathbb{P}^*_{\{x,y,\varphi,\hat{R}\}} = \frac{1}{S}\sum_{i=1}^{S} \delta_{\xi_{(i)}^{\{x,y,\varphi,\hat{R}\}}}$ *is the worst-case distribution, i.e.*

$$\inf_{\mathbb{P}\in \mathcal{F}_W(\theta)} \mathbb{E}_{\mathbb{P}}\left[h(x,\varphi,\widehat{R},\xi)\right] = \mathbb{E}_{\mathbb{P}^*_{\{x,y,\varphi,\widehat{R}\}}}\left[h(x,\varphi,\widehat{R},\xi)\right]$$

In other words, we can recover, for a given portfolio decisions, the worst-case distribution of $\xi^{\{x,y,\varphi,\widehat{R}\}}$, which is a Dirac distribution supported on a single worst-case scenario. In particular, in accordance with the imposed independence assumption, the worst-case scenarios for uncertain rate of return can be obtained separately.

## 4 Problem decomposition and optimality conditions:

The mixed-integer planning problem can be divided into two steps of optimization, the first step determining the portfolio composition and the second step determining the investment weights. In this way, the problem can be decomposed into a finite number of quadratic programming subproblems for parallel computation. For example, the decision maker selects $k$ stocks from $n$ stocks to form a portfolio and will have $C_N^K$ options, and each option is a quadratic programming subproblem for determining the weights of k stocks:

$$\underset{x,\lambda,\nu}{maximize} -\frac{1}{s}\sum_{i=1}^{s} \nu_i - \lambda\theta - \frac{1}{2}\mathcal{A}\, x^T\widehat{\Sigma}' x \qquad (g)$$

$$\begin{aligned} s.t. \quad & \sum_{j=1}^{k} x_j = 1 \\ & \nu_i + (1+\varphi)\widehat{\xi}_i^T x - \varphi\widehat{R} \geq 0, \forall i \in [S] \\ & \nu_i - \widehat{\xi}_i^T x \geq 0, \forall i \in [S] \\ & -\frac{\lambda}{1+\varphi} \leq x_j \leq \frac{\lambda}{1+\varphi}, \forall j \in [K] \\ & \lambda \geq 0 \\ & x_j \in [0,1], \forall j \in [K] \end{aligned}$$

Where $\widehat{\Sigma}'$ is the covariance matrix of the selected $k$ stocks, $\widehat{\Sigma}' \succcurlyeq 0$.

Introducing the decision variable $z = [x\ \nu\ \lambda]^T = [x_1,...,x_K, \nu_1,...,\nu_S, \lambda]^T$, and make $I_1^T z = \sum_{j=1}^{k} x_j$, $I_2^T z = \sum_{i=1}^{s} \nu_i$, $I_3^T z = \lambda$, which means:

$$I_1^T = \left[\underbrace{1,1,...,1}_{k\,1's}\underbrace{0,0,...,0}_{s+1\,0's}\right]^T, I_2^T = \left[\underbrace{0,0,...,0}_{k\,0's}\underbrace{1,1,...,1}_{s\,1's},0\right]^T, I_3^T = \left[\underbrace{0,0,...,0}_{k+s\,0's},1\right]^T$$

Then, the subproblems can be reformulated as:

$$\underset{z}{maximize} - \left[\frac{I_2}{S} + \theta I_3\right]^T z - \frac{1}{2}\mathcal{A}\, z^T Q z \qquad (h)$$

$$\begin{aligned} s.t \quad & G_i(z) \leq 0,\ i \in I = \{1,2,...,l\} \\ & H_i(z) = 0,\ i \in \varepsilon = \{1,2,...,m\} \end{aligned}$$

Where, $G(z) \leq 0$ is the inequality constraint, and $H(z) = 0$ is the equation constraint. $Q = \begin{bmatrix} \widehat{\Sigma}' & O_{k\times(s+1)} \\ O_{(s+1)\times k} & 0 \end{bmatrix}$ is $(k+s+1) - \text{order}$ split matrix, and $Q \succcurlyeq 0$ is semi-positive definite matrix. The $Karush - Kuhn - Tuck\ (kkt)$ conditions of the problem can be expressed as:

$$\begin{cases} Qz^* + \frac{I_2}{S} + \theta I_3 - \sum_{i=1}^{m} \mu_i^* \nabla H_i(z^*) - \sum_{i=1}^{l} \omega_i^* \nabla G_i(z^*) = 0 \\ H_i(z^*) = 0,\ i \in \varepsilon \\ G_i(z^*) \leq 0, \omega_i^* \geq 0, \omega_i^*\, G_i(z^*) = 0 \quad i \in I \end{cases}$$

Where $\mu$ and $w$ are auxiliary variables.

On the above convex quadratic programming problem, the common methods are: the Interior Point Method, the Active Set Method, Gradient Projection Method, etc. However, in real-life applications, if the accuracy of the solution is not required, we prefer to sacrifice some accuracy for faster solution efficiency, so the heuristic algorithm is widely used for solving large-scale NP-hard problems. In the later sections we design a hybrid algorithm to reduce the time cost of solving while maintaining a certain level of accuracy.

## 5 Hybrid Algorithm for DRP:

5.1 Tabu-Search algorithm：

We assume that the decision maker obtains maximum objective at the optimal solutions $x^* \in [0,1]$, $y^* \in \{0,1\}$ . Then, the problem $(f)$ is transformed into a simple linear programming.

$$\underset{\lambda,\nu}{maximize} -\frac{1}{S}\sum_{i=1}^{S} \nu_i - \lambda\theta - C$$

$$\begin{aligned} s.t. \quad & \nu_i \geq c_i^1, \forall i \in [S] \\ & \nu_i \geq c_i^2, \forall i \in [S] \\ & \|x^*\|_\infty \leq \frac{\lambda}{1+\varphi} \\ & \lambda \geq 0 \end{aligned}$$

Where $C = \frac{1}{2}\mathcal{A}\, x^{*T}\hat{\Sigma}' x^*$, and $\hat{\Sigma}'$ is the covariance matrix of the K assets determined by the decision variable $y^*$, $c_i^1 = \varphi\hat{R} - (1+\varphi)\hat{\xi}_i^{'T} x^*$, $c_i^2 = \hat{\xi}_i^{'T} x^*$.

Obviously, in order to make the objective function take the maximum value, we should have $\lambda = (1+\varphi) \cdot \|x^*\|_\infty$, and

$$\text{For } \forall i \in [S],\ v_i = \begin{cases} \nu_i = \varphi\hat{R} - (1+\varphi)\hat{\xi}_i^{'T} x^*, & c_i^1 \geq c_i^2 \\ \nu_i = \hat{\xi}_i^{'T} x^*, & c_i^1 < c_i^2 \end{cases}$$

Then, the max objective value is

$$f(x^*,y^*) = -\frac{1}{S}\sum_{i}^{S} \max\{c_i^1, c_i^2\} - \theta(1+\varphi) \cdot \|x^*\|_\infty - C$$

Starting from an initial feasible solution $\check{x}_0^T \in [0,1]$, $\check{y}_0^T \in \{0,1\}$, a series of specific search directions are chosen as trials to achieve an iterative optimization that allows a specific objective function value. To avoid getting trapped in a local optimum, a flexible "memory" technique is used in the Tabu search to record and select the optimization process that has been performed and to guide the next search direction.

The process about the Tabu-Search for DRP is shown in Tab.1:

**Tab.1 Tabu-Search Algorithm**

| **Tabu-Search** |
| --- |
| (**Initialization**) |
| • Initial test solutions: $\check{x}_0^T \in [0,1]$, $\check{y}_0^T \in \{0,1\}$; bestCandidate: X ← $\check{x}_0^T$; Y ← $\check{y}_0^T$. |
| • Initial objective value: bestObj ←− inf |
| • Setting **TabuList** ← [ ];$TabuList.add$(X,Y) |

```
• Setting Ignoring Rule
(Iteration)
  While (not stopping conditions ())
    nbX ← genNeighbor(X,n); nbY ← genNeighbor(Y,n)
    nbObj ← NeighborObj(nbX,nbY)
    For i = 1:n
      If (nbObj(i) > bestObj) && (nbX(i),nbY(i) ∉ TabuList)
        bestObj ← nbObj(i);
                X ← nbX(i), Y ← nbY(i)
                TabuList.add(X,Y)
      Else
        Continue
    End
    For i = 1:length(TabuList)
      If (satisfied Ignoring Rule)
                TabuList.release(TabuList(i))
      Else
        Continue
    End
  End
• Return bestObj, bestCandidate: X,Y
```

Italicized words are functions:

- Func1>>*genNeighbor*(■,n)*: Based on the solution ■ generate n neighborhood solutions. The set of neighborhood solutions of Y is generated by the single point exchange, and the determination of the exchange points obeys some distribution.*
- Func2>> *NeighborObj : Under the given solution Y conditions, solve the optimal function values sequentially by neighborhood solutions.* $f(x,y) = -\frac{1}{S}\sum_i^S max\{c_i^1,c_i^2\} - \theta(1+\varphi) \cdot \|x\|_\infty - C$
- Func3>>*TabuList.add: Add the candidate solutions to the TabuList.*
- Func4>>*TabuList.release: Release the candidate solutions from the TabuList.*

In general, the Tabu-Search algorithm is difficult to guarantee the exact solution of the problem. Since it is by nature a search-on-the-fly, and although it can avoid repeated searches to some extent, the problem may easily fall into local solutions with iterations. Therefore, this algorithm is mainly suitable for the case where the accuracy of the problem solution is not required, and a rational termination criterion or iteration step should be determined. To improve the algorithm accuracy, we give a hybrid algorithm strategy

5.2 Hybrid Algorithm:

The idea of the hybrid algorithm is to take advantage of the Tabu-Search algorithm for the decision variable $y$, then determine the subproblem $(h)$ and construct the penalty function of the objective function for optimization.

First, the penalty function of the subproblem $(h)$ is constructed under the solution $y^T \in \{0,1\}$. By introducing the relaxation variable $g_j$,$i = 1,2,...,l$, the problem $(h)$ is equivalently transformed into:

$$
\begin{aligned}
&\underset{z}{maximize} - \left[\frac{I_2}{S} + \theta I_3\right]^T z - \frac{1}{2}\mathcal{A}\, z^T Q z && (I)\\
&s.t \qquad G_i(z) - g_i = 0,\ i \in I = \{1,2,...,l\}\\
&\qquad\quad\ H_i(z) = 0,\ i \in \varepsilon = \{1,2,...,m\}\\
&\qquad\quad\ g_i \geq 0,\ i \in I = \{1,2,...,l\}
\end{aligned}
$$

Then, construct the mixed objective expansion function for the equivalence problem $(I)$.

$$\psi(z,g,\tau) = -\left[\frac{I_2}{S} + \theta I_3\right]^T z - \frac{1}{2}\mathcal{A}\, z^T Q z + \frac{1}{2\tau}\sum_{i=1}^{m} H_i^2(z) + \frac{1}{2\tau}\sum_{i=1}^{l} [G_i^2(z) - g_i]^2 - \tau\sum_{i=1}^{l} ln g_i$$

On this basis we can establish the corresponding algorithm for solving the unconstrained optimization problem. Any $x$, $g$ $(g > 0)$ can be used as a suitable initial point to start the iterative algorithm, and the flow of the hybrid algorithm is shown in Tab.2.

**Tab.2 Hybrid Algorithm**

**Hybrid Algorithm**

(**Initialization**)
- Initial test solutions: $\check{y}_0^T \in \{0,1\}$; bestCandidate: $Y \leftarrow \check{y}_0^T$.
- Initial objective value: bestObj ←− inf
- Setting **TabuList** ← [ ];*TabuList.add*(Y)
- Setting **Ignoring Rule**

(**Iteration**)
**While** (not stopping conditions ())
  nbY ← *genNeighbor*(Y,n) ,
  **For** i = 1:n
    Initial penalty factor $\tau^{(0)} > 0$; Allowable error $\varepsilon > 0$
    Initial Point $z^{(k-1)} \longleftarrow z^0, g^{(k-1)} \leftarrow g^0$   k=1
    Constructing subproblem penalty functions:
    $\psi(z,g,\tau) = -\left[\frac{I_2}{S} + \theta I_3\right]^T z - \frac{1}{2}\mathcal{A}\, z^T Q z + \frac{1}{2\tau}\sum_{i=1}^{m} H_i^2(z) + \frac{1}{2\tau}\sum_{i=1}^{l}\left[G_i^2(z) - g_i\right]^2 - \tau\sum_{i=1}^{l} \ln g_i$
    **While** $\left(\left\| z^*\left(\tau^{(k)}\right) - z^*\left(\tau^{(k-1)}\right)\right\| \le \varepsilon\right)$
      min $\psi(z,g,\tau)$ for $z^*\left(\tau^{(k)}\right)$; $k = k+1$; $\tau^{(k+1)} = C\tau^{(k)}$
    nbObj(i) = $\psi(z^*, g^*, \tau^*)$
  **End**
  nbObj ← $z^*$n,nbY)
  **For** i = 1:n
    **If** (nbObj(i) > bestObj) && (nbY(i) ∉ **TabuList**)
      bestObj ← nbObj(i);  Y ← nbY(i); *TabuList.add*(Y)
    **Else**
      Continue
  **End**
  **For** i = 1:*length*(**TabuList**)
    **If** (satisfied **Ignoring Rule**)
        *TabuList.release*(**TabuList**(i))
    **Else**
      Continue
  **End**
**End**
- Return bestObj, bestCandidate:Y

The hybrid algorithm is very efficient in solving DRP and is also applicable to other mixed integer convex optimization problems.

**Tab.3 Comparison of computational efficiency**

| ***Ambiguity Size*** | ***Solver* by Cplex** | | ***TabuSearch* (2000 *iter*)** | | ***Hybrid Algorithm* (300 *iter*)** | |
|---|---|---|---|---|---|---|
| (Init $\theta = 0.001$) | *ObjValue* | *Time (s)* | *ObjValue* | *Time (s)* | *ObjValue* | *Time (s)* |
| **1** × | -0.009864 | 675.157903 | -0.010856 | 301.163248 | -0.010201 | 420.809468 |
| **2** × | -0.010206 | 1500.731880 | -0.011144 | 302.139409 | -0.010585 | 413.470442 |
| **3** × | -0.010499 | 1757.639750 | -0.011782 | 322.638795 | -0.011044 | 411.943905 |
| **4** × | -0.010760 | 2407.372070 | -0.012411 | 299.254264 | -0.011712 | 430.822711 |
| **5** × | -0.011016 | 3536.957276 | -0.012669 | 331.600296 | -0.011905 | 430.710590 |

We selected the daily returns of NDX 100's constituent stocks from 2010.01.01-2020.01.01 as experimental subjects to test the algorithm performance (see Tab.3).

It has been shown that with the expansion of the ambiguity set, the hybrid algorithm can significantly reduce the time cost of the solution, while maintaining a high solution accuracy compared to the classical Tabu-Search algorithm.

## 6 Sensitivity Analysis of Emotional Parameters:

The correlation aversion coefficient ($\varphi \geq 0, A \in R, \theta \geq 0, \widehat{R} \in \mathfrak{R}$) will affect the utility of investors. For this purpose, we conducted experiments on the robustness of the model, and the sensitivity analysis for the target utility is shown in Fig.1-6

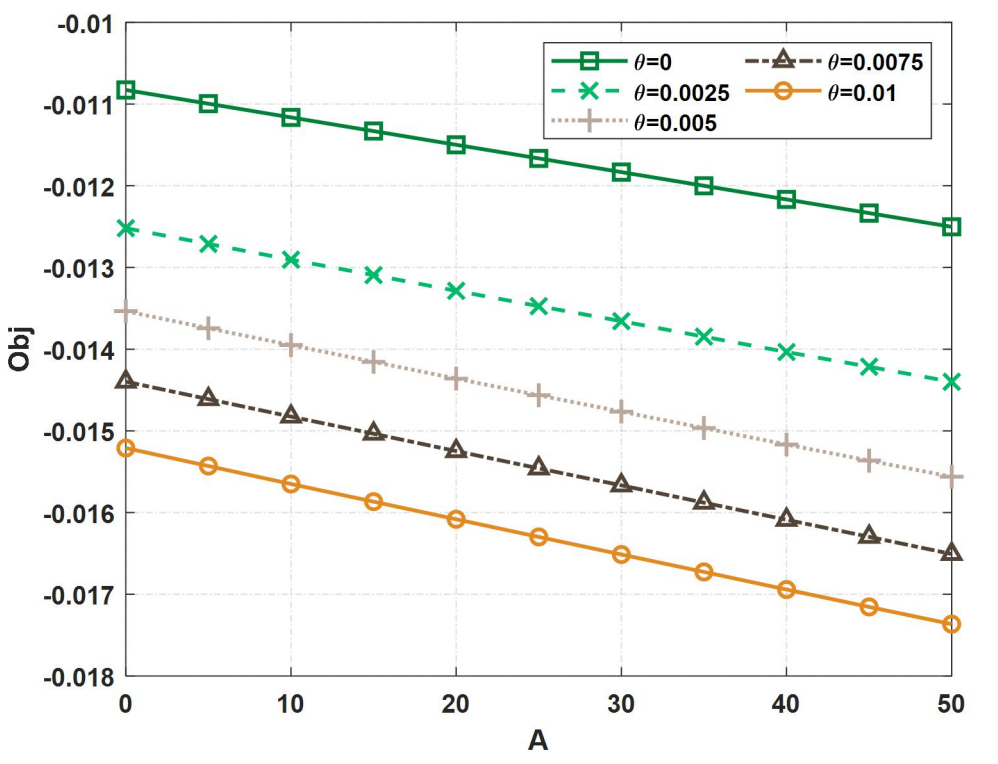

**Fig.1-1 Ambiguity Aversion & Risk Aversion**

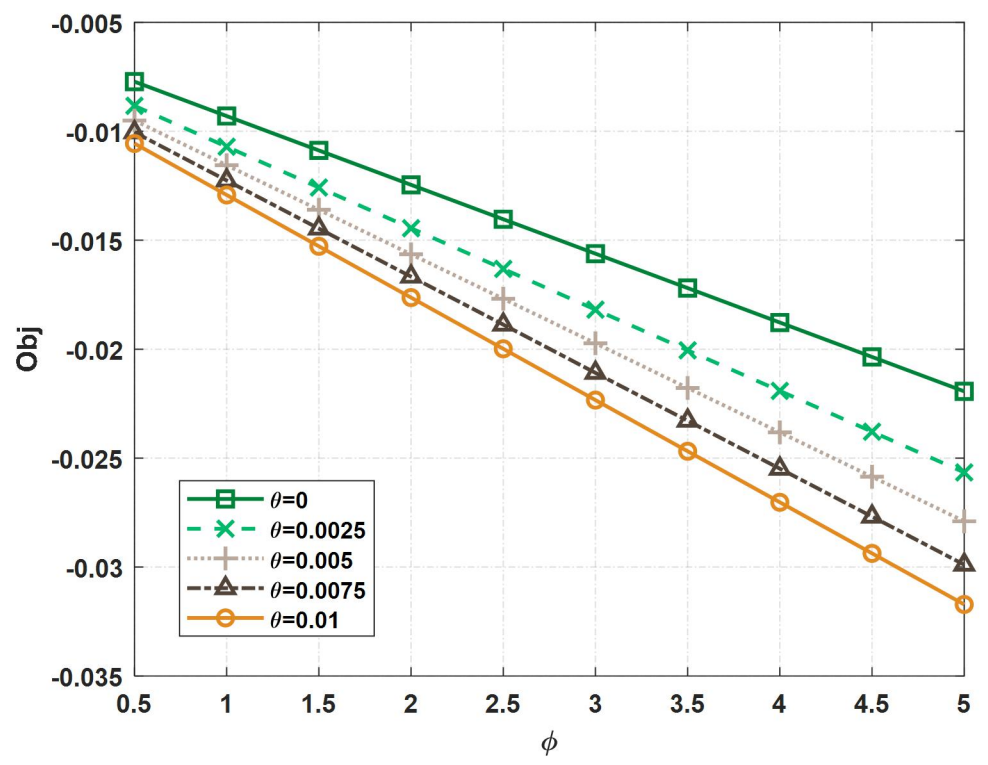

**Fig.1-2 Ambiguity Aversion & Loss Aversion on**

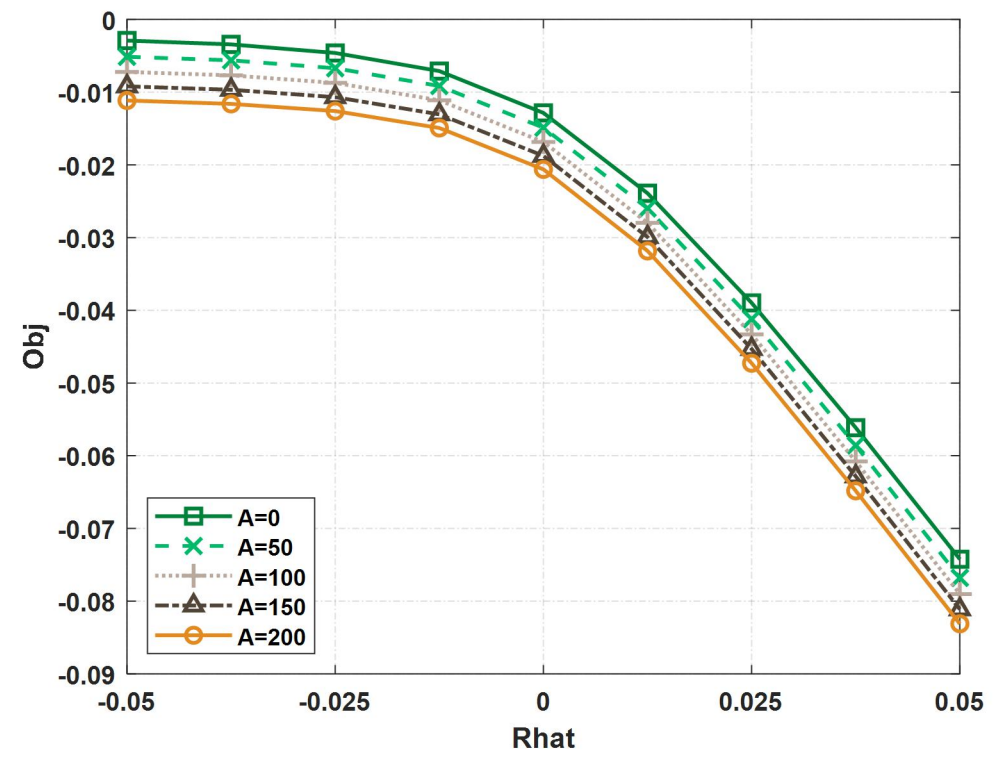

**Fig.1-3 Risk Aversion & Loss Reference Points**

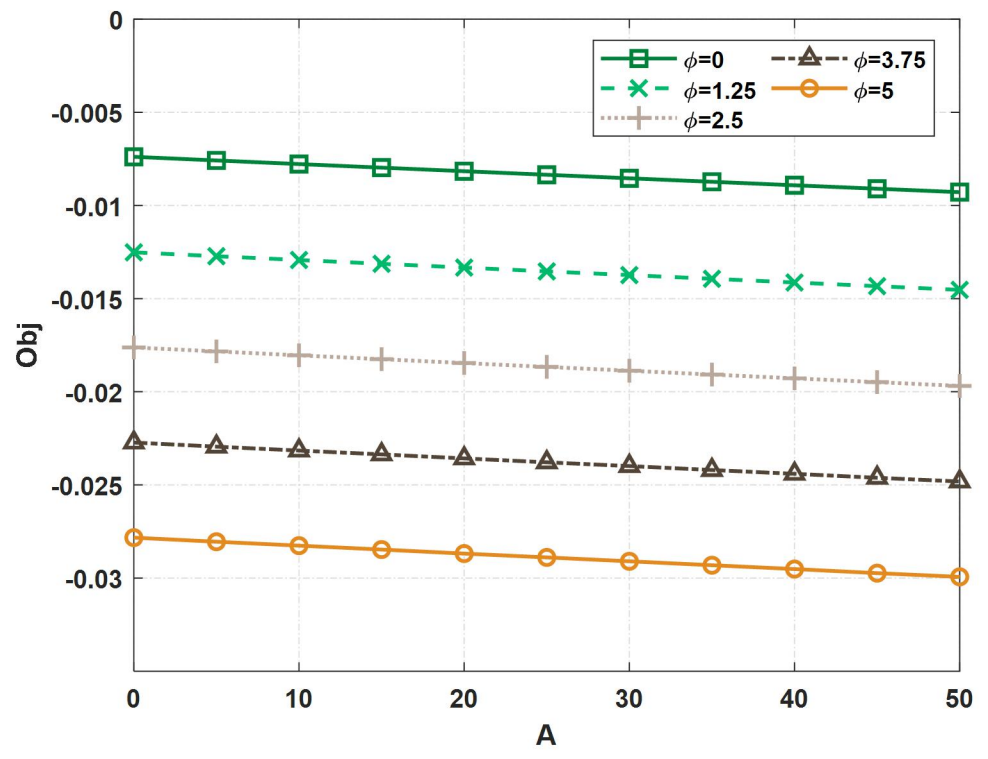

**Fig.1-4 Loss Aversion & Risk Aversion**

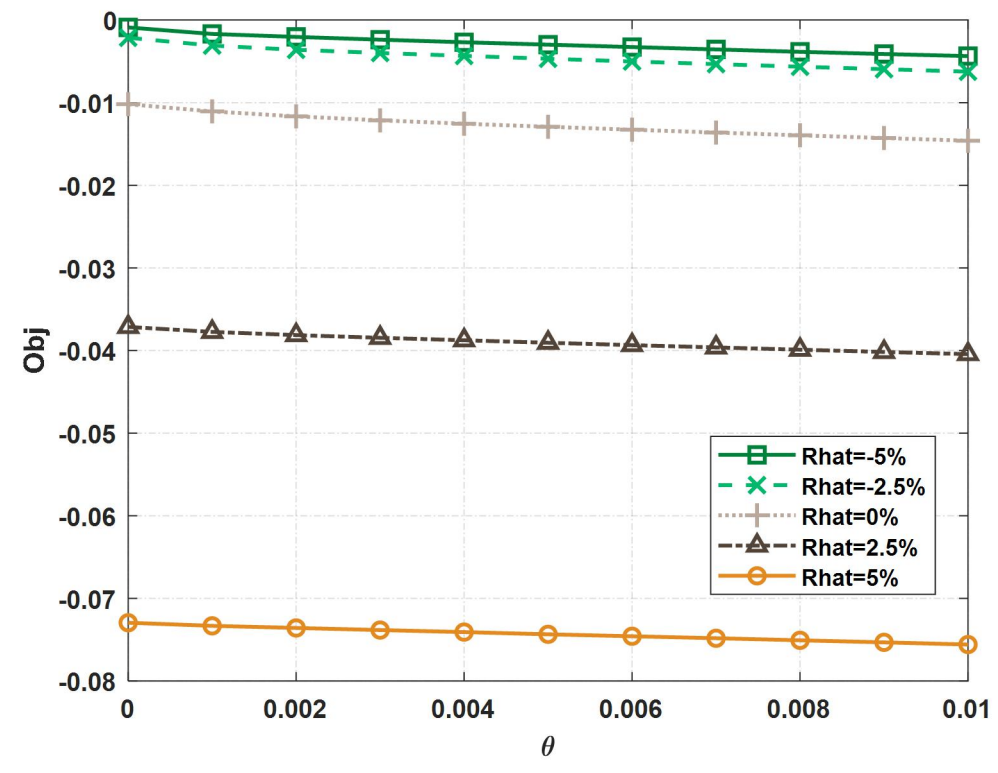

**Fig.1-5 Ambiguity Aversion & Loss Reference Points**

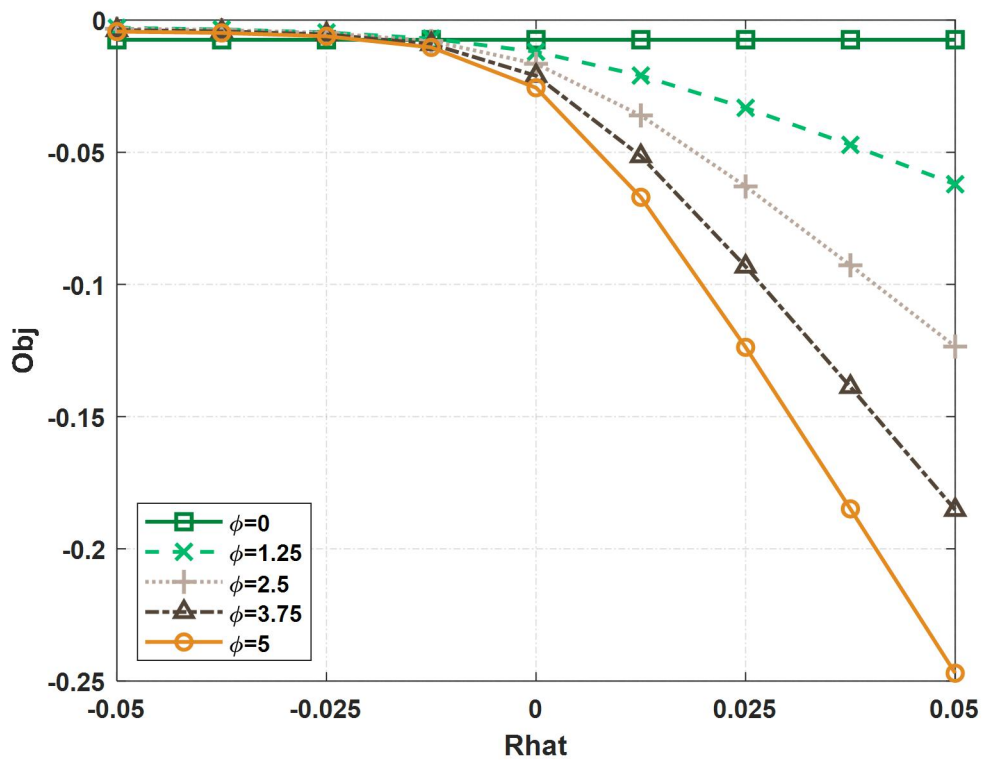

**Fig.1-6 Loss Aversion & Loss Reference Points**

The experimental result shows that the investors' aversion reduces the utility, and the effect of the aversion coefficient is linear with little crossover effect across different aversions. However, the investors' choice of loss reference points exacerbates the degree of impact of aversion on utility. That is, the higher the investor's expectations of the given investment, the more significant the effect of aversion.

Loss aversion reflects the fact that people do not have the same risk preferences; when it comes to gains, people are risk averse; when it comes to losses, people behave as risk seekers. We introduce three types of dynamic investors in a continuous investment process with out-of-sample forecasts under rolling time windows and consider the different changes in loss aversion of investors when they suffer losses to analyze the impact of different loss aversion characteristics of investors on portfolio utilities.

The dynamically changing loss aversion coefficients and reference points for the three types of investors are expressed as:

• **Type 0:** Investors' loss aversion is static. Their loss aversion coefficients $\varphi$ and reference points $\widehat{R}$ do not change due to changes in wealth.

• **Type 1:** Their dynamically changing loss aversion coefficients and reference points are expressed as:

$$\lambda_t = \begin{cases} \lambda_0, & W_t \geq W_{t-1} \\ \lambda_0 + \left(\frac{W_{t-1}}{W_t} - 1\right), & W_t < W_{t-1} \end{cases} \qquad \widehat{R}_t = \begin{cases} \widehat{R}_0, & W_t \geq W_{t-1} \\ \frac{W_{t-1}}{W_t}\widehat{R}_0, & W_t < W_{t-1} \end{cases}$$

Where, $W_t$ is the wealth value at time $t$, $\lambda_0$ and $\widehat{R}_0$ are the investor's initial loss aversion coefficient and reference point, respectively, and the loss aversion coefficient $\lambda_t$ at time $t$ satisfies $\lambda_t > \lambda_0 > 0$ and the reference point $\widehat{R}_0$ at time $t$ satisfies $\widehat{R}_t > \widehat{R}_0 > 0$. The investors have increased loss aversion coefficients and reference points When the moment is suffering losses. They are relatively conservative in their investment strategy and may be reluctant to sell losing assets due to fear of further losses.

• **Type 2:** Their dynamically changing loss aversion coefficients and reference points are expressed as:

$$\lambda_t = \begin{cases} \lambda_0 + \left(\frac{W_t}{W_{t-1}} - 1\right), & W_t \geq W_{t-1} \\ \lambda_0, & W_t < W_{t-1} \end{cases} \qquad \widehat{R}_t = \begin{cases} \frac{W_t}{W_{t-1}}\widehat{R}_0, & W_t \geq W_{t-1} \\ \widehat{R}_0, & W_t < W_{t-1} \end{cases}$$

Similarly, the loss aversion coefficient $\lambda_t$ at time $t$ satisfies $\lambda_t > \lambda_0 > 0$ and the reference point $\widehat{R}_0$ at time $t$ satisfies $\widehat{R}_t > \widehat{R}_0 > 0$. For this group of investors, losses that have already occurred do not exacerbate their emotional aversion, but when wealth is growing, they fear those unseen losses in the future, which manifest themselves in a rush to sell assets that are rising.

Sensitivity analysis of the dynamic loss aversion coefficient is based on the daily returns of the NDX100's stocks. Dynamic investors select stocks based on historical data, and during the holding period, decision makers' loss aversion accumulates (manifested as a reduction in utility) and leads to frequent trading and changes in the asset weights of the positions. Alternatively, this is a continuous process of position adjustment during the holding period.

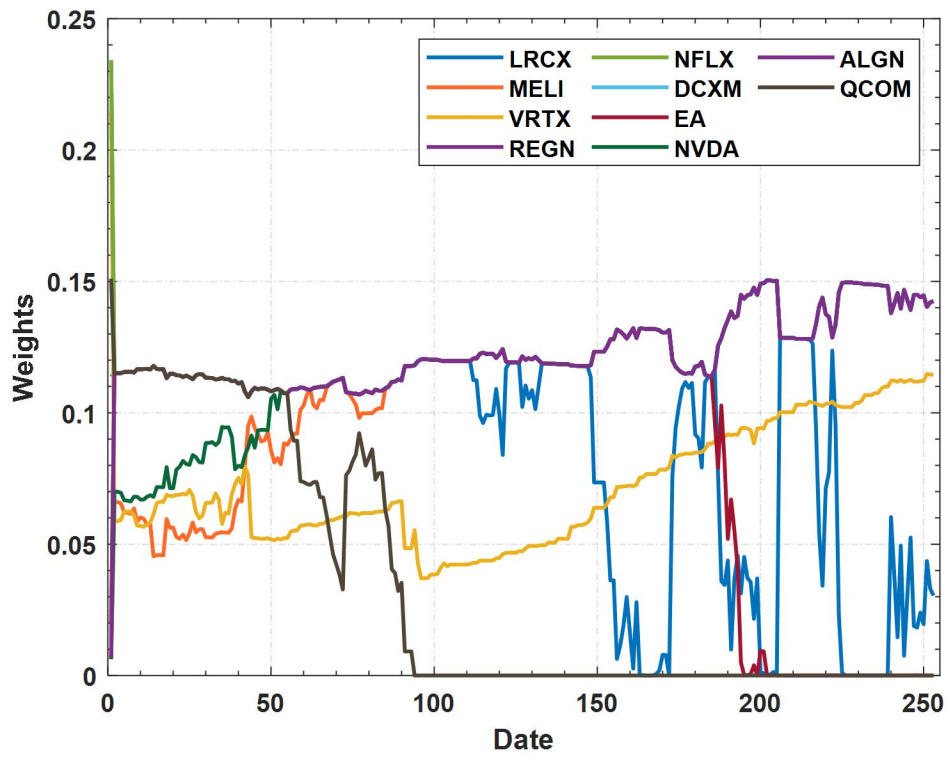

**Fig.2-1 Type-1 Dynamic Investors' Decision Weights**

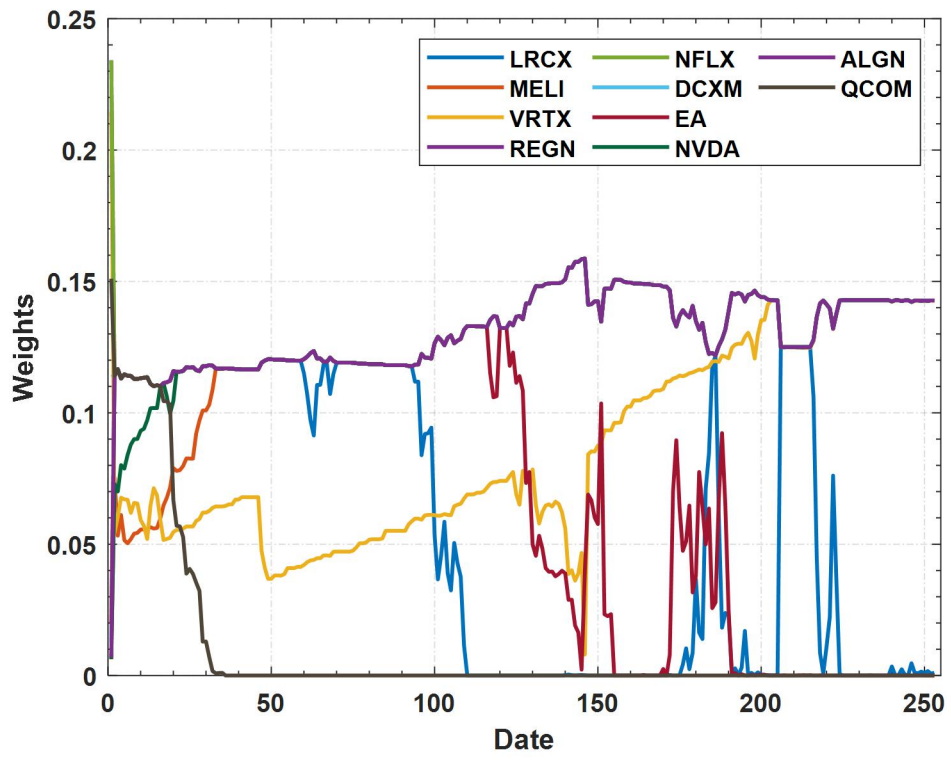

**Fig.2-2 Type-2 Dynamic Investors' Decision Weights**

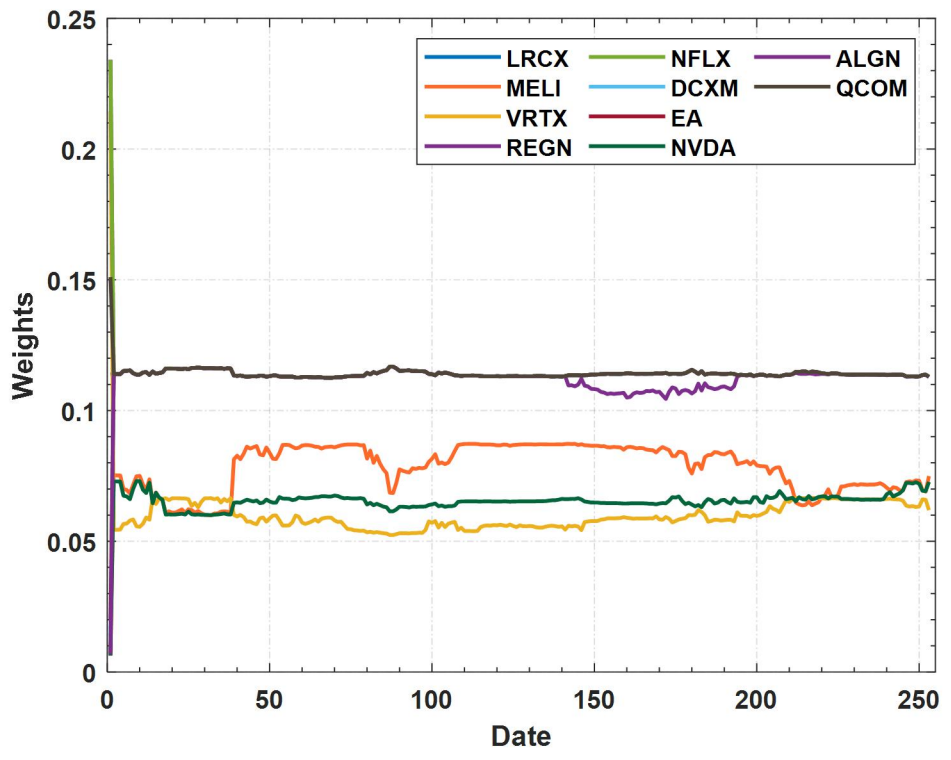

**Fig.2-3 Type-0 Dynamic Investors' Decision Weights**

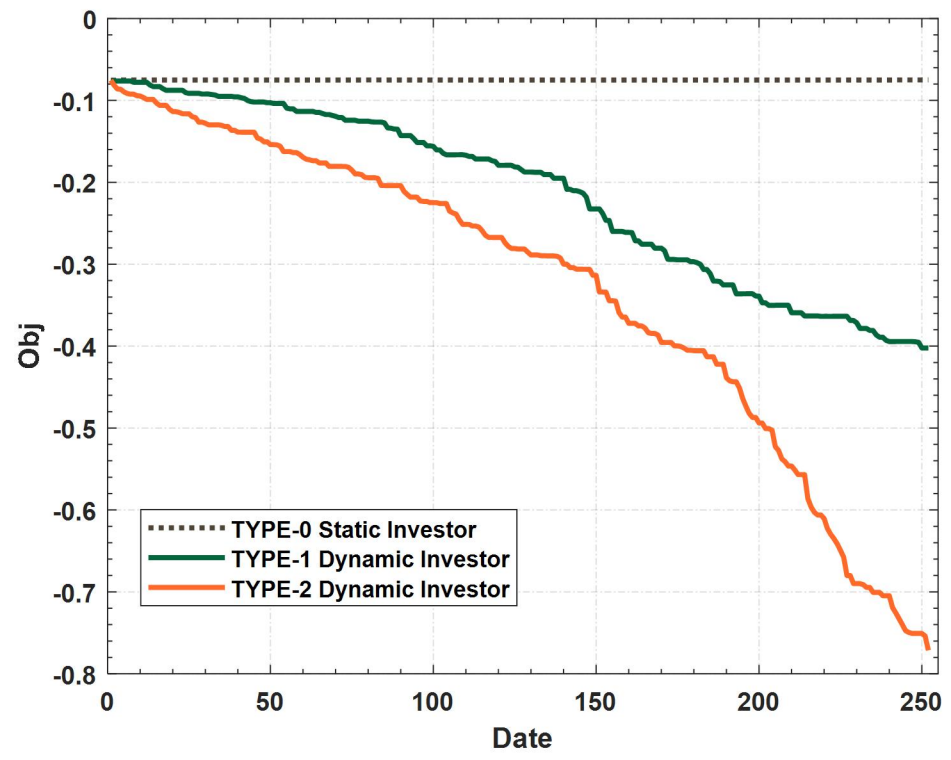

**Fig.2-4 Utility Objectives of Three Types of Dynamic Investors**

Dynamic loss averse investors adjust their positions more frequently and to a greater extent than static investors. And as aversion accumulates, Type-2 investors' utility is cut significantly faster than type-1 investors. This also indicates to some extent that loss aversion utility is asymmetric

## 7 Empirical Research:

### 7.1 Simulation of Portfolio Selection

Based on monthly return data for NDX100's stocks from Jan.2015-Jan.2020, the dataset is divided into five period windows for out-of-sample portfolio return forecasting. Setting investor's aversion characteristics: loss aversion coefficient $\varphi = 1.5$ , daily profit/loss reference point $\widehat{R} = 0.1\%$, risk aversion coefficient $\mathcal{A} = 1.5$, ambiguity aversion coefficient $\theta = 0.003$, the maximum iteration of the algorithm is set to 300. We examined wealth growth and returns over the back-test period separately for portfolios of different sizes (from 5 to 80 stocks), see fig.3.

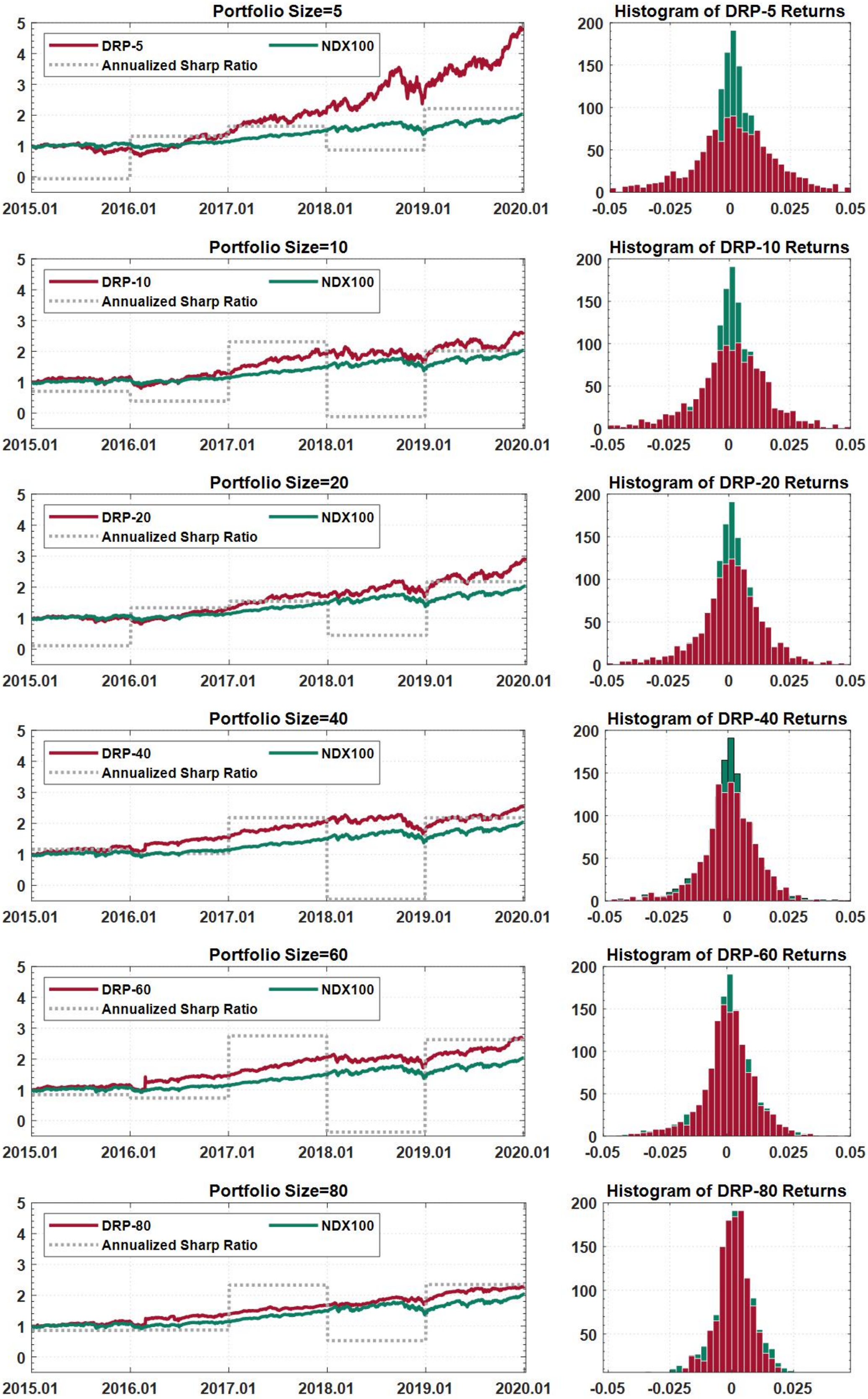

**Fig.3 Risk diversification becomes higher when the size of the portfolio expands**

Of course, since the DRP problem is a non-convex quadratic programming problem, and the optimal solution of the hybrid solution algorithm is not globally optimal, we compare some common portfolio construction strategies in order to reflect the advantages of DRP in asset allocation.

7.2 Comparison of Asset Allocation Model

In some unexpected events the market fluctuates sharply, accompanied by the possibility of large declines, and it is difficult to hedge these losses based on historical data performance, for this reason the empirical study compares the performance of the DRP model relative to common asset allocation strategies since the outbreak of the COVID-19.

Data source:

The relevant dataset for asset allocation efficiency comparisons is obtained from the Wind database. We select 474 stocks with longevity in the SP500 index. The initial decision is made based on the daily return from 2015.01.01 to 2019.9.30, with a back-test point of 2019.9.30 (the wealth of the back-test point is set to one), and thereafter, the decision maker reallocates the asset weights once a year.

Asset allocation strategies during the COVID-19:

The process by which investors employ some common asset allocation strategies (Tab.5) to achieve wealth during the global financial market shakeout caused by the Corona Virus Disease 2019. Wealth changes and the SP500 index are plotted in Fig. 4.

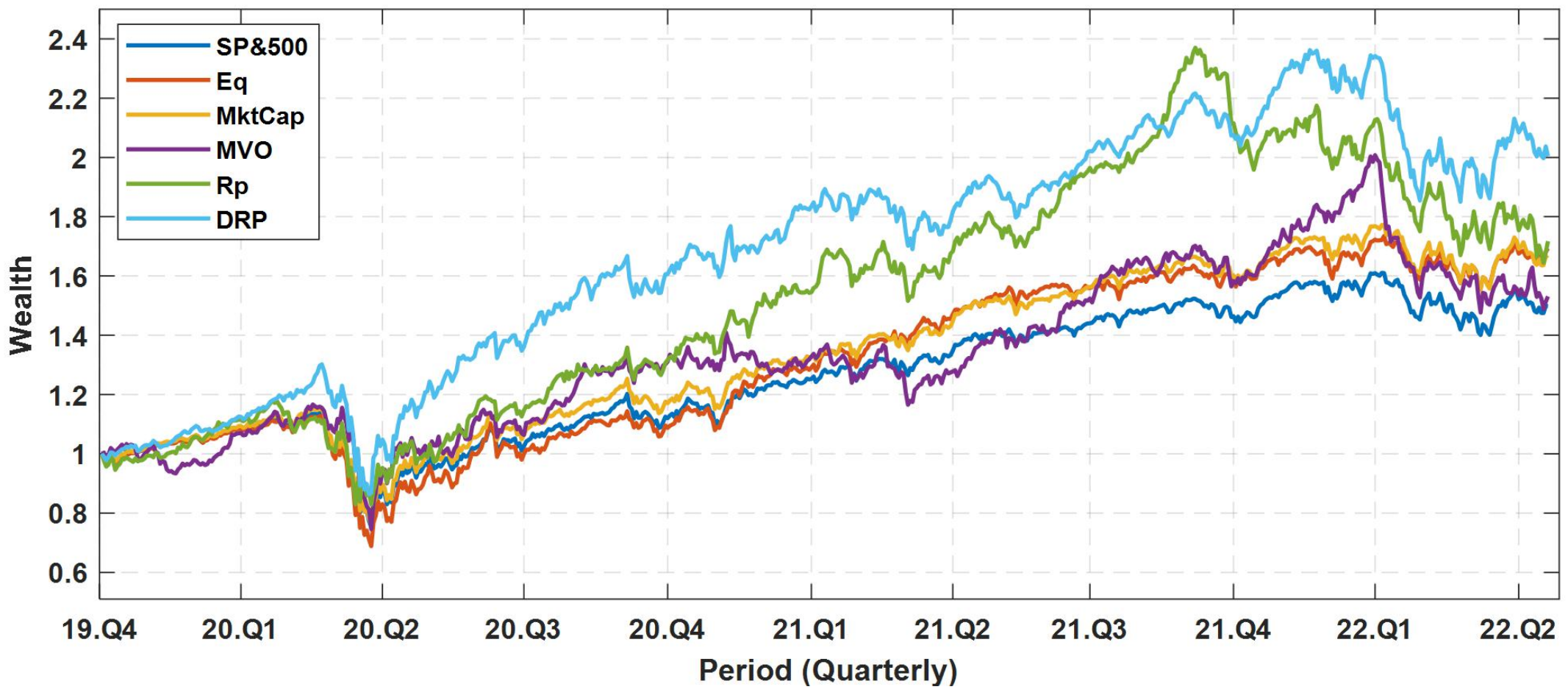

**Fig.4 Asset Allocation Model Wealth Comparison**

The global outbreak and spread of the Corona Virus Disease (COVID-19) had a severe impact on the financial markets. There was a very sharp correction in global stock markets, with U.S. stocks melting down several times in a few weeks, volatility in financial markets spiking significantly, and the market also had a rare phenomenon of risky and safe-haven assets falling in tandem. The main reason for this was the unanticipated spread of the new crown pneumonia epidemic overseas, which caused a change in investors' risk appetite and expectations and triggered a liquidity run on financial markets. Graphically, since the outbreak of the epidemic COVID-19 in 2019.10 to date, financial markets have grown slowly and accompanied by declines. The DRP strategy significantly outperformed common asset allocation strategies, as demonstrated by the fact that at the end of the first quarter of 2020 (20.Q1), when U.S. stocks experienced a meltdown, our strategy managed downside losses as much as possible and outperformed the index and other asset allocation strategies in the subsequent rebound phase.

**Tab.4 Asset Allocation Model**

| Strategy | Definition | Advantages | Disadvantages |
|---|---|---|---|
| **Equal Weighting** | Assign equal weights to each constituent asset; | (1). Avoid over-representation of some assets to fully diversify risk;<br>(2). Can passively sell assets with large short-term increases and buy assets with large short-term decreases; | (1). Increased liquidity, not suitable for large-scale portfolio, when the portfolio size is relatively large, the smallest size of assets, it will limit the liquidity and size of this portfolio;<br>(2). Larger elasticity, amplify market fluctuations up and down; |
| **Market-Value Weighting** | Assume that investors determine the investment weights in proportion to market capitalization | (1). No need for active management, market capitalization weighted portfolio automatically completes individual stock allocations when stock prices fluctuate;<br>(2). Ensures that the majority of the portfolio is allocated to liquid assets due to the liquidity of companies with large market capitalization, reducing shock costs.<br>(3). Higher Sharpe ratio. | (1). Vulnerable to large cap stocks, assigning higher weights to assets with pricing errors and potentially suffering more significant declines when stock prices are not efficient |
| **Risk Parity** | Risk Parity is an asset allocation concept that assigns equal risk weights to different assets in a portfolio | (1). Reduces the portfolio risk of traditional allocation methods by making the risk of each asset class more evenly distributed in the portfolio | (1). Lower returns at certain times of the year, which do not satisfy some investors' preference for higher risk and higher returns |
| **Markowitz MVO** | The portfolio in which the investor obtains the maximum expected utility among all possible portfolios based on mathematical and statistical methods. Determine the portfolio that maximizes return for a given level of risk, or minimizes risk for a given level of return | (1). Takes into account the risky asset preferences are well compatible with Sharpe theory.<br>(2). Modeling long-term returns can solve the problem of mismatch between input maturity and portfolio selection maturity. | (1). MVO assumes that asset returns are normally distributed, but actual asset returns exhibit significant non-normal characteristics;<br>(2). MVO does not give a method for generating optimal portfolios;<br>(3). Estimation bias can lead to non-effective portfolios;<br>(4). MVO is a single investment period model, which does not connect investors well across investment horizons and may lead to non-optimal returns as investment horizons become longer<br>(5). No consideration of liquidity |

In order to be able to compare the different strategies more visually, we also compare the performance of different portfolios over the period, including: annualized expected return, Sharpe ratio, maximum retracement, Jensen alpha, Traynor ratio, information ratio, etc. The five-year annualized performance and averages are shown in Tab.5, where the underlined values indicate that the DRP model significantly outperforms market benchmarks and general asset allocation strategies in terms of exposure to systematic risk and excess returns.

**Tab.5 Asset Allocation Model Comparison**

| | SP&500 | Eq | MV | MVO | Rp | DRP |
|---|---|---|---|---|---|---|
| **Annualized Rate of Return** | 0.1883 | 0.2420 | 0.2298 | 0.2122 | 0.2548 | **0.3138** |
| **Annualized Standard Deviation** | 0.2431 | 0.2649 | 0.2474 | 0.2993 | 0.2918 | 0.2859 |
| **Sharpe Ratio** | 0.7197 | 0.8632 | 0.8751 | 0.6644 | 0.8277 | **1.0509** |
| **Drawdown** | 38.59% | 43.86% | 39.38% | 53.11% | 72.62% | 51.17% |
| **Beta** | 1.0000 | 1.0399 | 1.0151 | 0.9179 | 0.9317 | **1.0539** |
| **Jensen's alpha** | 0.0000 | 0.0468 | 0.0389 | 0.0383 | 0.0785 | **0.1161** |
| **Treynor Ratio** | 0.1749 | 0.2199 | 0.2133 | 0.2166 | 0.2592 | **0.2851** |
| **Information Ratio** | -- | 0.6721 | 2.2805 | 0.1192 | 0.3603 | 0.9835 |

(*The market yield is referenced to the SP&500 Index and the risk-free rate is referenced to the US 10-year T-NOTES.)

## 8 Concluding Remarks:

We provide a distributed robust optimization portfolio selection model that takes into account multiple aversion characteristics of investors. The distributed robust strategy reduces the sensitivity of traditional optimized portfolios to returns and also avoids traditional robust optimization that is too conservative and impractical. In practice, such a strategy can be an effective hedge against the systematic risk of a position when financial markets fluctuate sharply. Our work also includes a semi-heuristic hybrid algorithm design, which allows for controlled accuracy in model solving and significant time savings in applying to large-scale problems.

Some possible extensions include whether the covariance between assets can also be taken into account for ambiguity; how the model will perform in some markets containing mispriced assets; Considering bull and bear markets separately; Algorithm improvement.

## Reference：